\newcommand{\C}{{\bf C}}
\newcommand{\rlh}{\rightleftharpoons}
\newcommand{\dar}{\downarrow}
\newcommand{\al}{\alpha}
\newcommand{\La}{L_{\alpha}}

\documentclass[12pt]{amsart}

\begin{document}

\author{Victor~A.~Vassiliev}

\email{vva@mi.ras.ru}

\thanks{
Mathematics College of Independent Moscow University \& Steklov Mathematical
Institute Supported in part by the Russian Fund for Basic Investigations
(Project 95-01-00846a) and INTAS grant (Project 4373).}

\title{Stratified Picard-Lefschetz Theory with Twisted Coefficients}

\date{Revised version was published in 1997 and dedicated
to V.I.~Arnold on the occasion of his 60th anniversary}

\begin{abstract}
The monodromy action in the homology (generally with twisted coefficients) of
complements of stratified complex analytic varieties depending on parameters is
studied. For a wide class of local degenerations of such families (stratified
Morse singularities) local monodromy and variation operators are reduced to
similar operators acting in the transversal slices of corresponding strata.
These results imply a main part of generalized Picard--Lefschetz formulae of
\cite{Ph1} and (in the case of constant coefficients) similar reduction
theorems of \cite{V1}, \cite{V2}.
\end{abstract}

\maketitle

\section*{Introduction}

Let $M^n$ be a complex analytic manifold and $\Lambda_{\tau}$ a family of
subvarieties in $M^n$ depending analytically on the parameter $\tau$. For
almost all values of $\tau$ the corresponding pairs $(M^n,\Lambda_\tau)$ are
homeomorphic to one another, and the fundamental group of the set of such
generic $\tau$ acts on different homology groups related to such pairs. This
action is responsible for the ramification and qualitative analytic features of
all known functions given by integral representations (such as Fourier and
Radon transforms, fundamental solutions of hyperbolic and parabolic equations,
Newton--Coulomb potentials, hypergeometric functions, Feynman integrals etc.).
Explicit formulae expressing this action are called (generalized)
Picard--Lefschetz formulae, see e.g. \cite{AVG}, \cite{AVGL}, \cite{Ph1,Ph2},
\cite{giventh}, \cite{V1}.

The classical Picard--Lefschetz formula describes the case when $\Lambda_\tau$
is the family of irreducible hypersurfaces in $M^n$, non-singular for generic
$\tau$ and having a Morse singularity for $\tau$ running over a hypersurface in
the parameter space, the investigated loop is a small circle embracing this
hypersurface, and the homology group in question is $H_{n-1}(\Lambda_\tau)$.

Similar formulae for many other possible degenerations of the variety
$\Lambda_\tau$ and other homology groups were studied in \cite{Ph1} (see also
\cite{giventh}) and in \cite{V2}. In the present paper we consider the most
general situation: as in \cite{V1}, \cite{V2} we investigate the {\it
stratified Morse singularities} of $\Lambda_\tau$ (which include all
degenerations of \cite{Ph1} as special cases), and as in \cite{Ph1}, \cite{Ph2}
we consider the homology groups of $M^n \setminus \Lambda_\tau$ (with twisted
coefficients, and with both closed and compact supports), to which the study of
almost all other homology groups related with the pair $(M^n, \Lambda_\tau)$
can be reduced (except may be for the intersection homology of $\Lambda_\tau$
studied in the last part of \cite{V2}).
\medskip

The Arnold's ``complexification'' functor (see e.g. \cite{A}), establishing an
informal analogy between objects of ``real'' and ``complex'' worlds, maps the
Morse theory  into the Picard--Lefschetz theory; certainly, the analogue of the
stratified Morse theory of \cite{GM} should be a stratified Picard--Lefschetz
theory. On the other hand, recently it became clear that instead of the
homology groups of (sub)varieties, usually considered in all these theories, in
the ``complex'' situation it is natural to deal with the homology, generally
with coefficients in non-constant local systems, of their complements, cf.
\cite{Ph2}, \cite{gvz}, \cite{giventh}. Thus the matter of our paper seems to
be the most adequate ``complexification'' of that of \cite{GM}.
\medskip

{\bf Agreements.} In what follows we assume that $M^n = \C^n$ (which is not
restricting because all our considerations are local), and families
$\Lambda_\tau$ (which appear later as $A \cup X_t$) depend on one parameter. We
consider only homology groups reduced modulo a point in the case of absolute
homology and modulo the fundamental cycle in the case of relative homology of a
complex analytic variety. We often use a short notation of type $H_*(X,Y)$
instead of a more rigorous $H_*(X, Y \cap X)$ or $H_*(X \cup Y, Y)$. The sign
$\Box$ denotes the end or absence of a proof.

\section{Main characters and stating the problem}

Let $A$ be a complex analytic subvariety in $\C^n$ with a fixed analytic
Whitney stratification (see e.g. \cite{GM}, \cite{V1}), let $\sigma \subset A$
be a stratum of dimension $k$, $a$ a point of $\sigma$, and $B \equiv
B_\varepsilon$ a small closed disc in $\C^n$ centred at $a$, so that $B$ has
nonempty intersections only with the strata adjoining $\sigma$ (or $\sigma$
itself) and $\partial B$ is transversal to the stratification. In particular
(see \cite{Che}, \cite{GM}) the induced partition of the pair $(B, \partial B)$
is again a Whitney stratification of $B$.

Let $f: (B,a) \to (\C,0)$ be a holomorphic function such that $df \ne 0$ in
$B$.
\medskip

{\bf Definition 1} (cf. \cite{GM}). The function $f$ has a {\it Morse
singularity} on $A$ at the point $a$ if its restriction on the manifold
$\sigma$ has a Morse singularity at $a$, and for any stratum $\tau \ne \sigma$
and any sequence of points $b_i \in \tau$ converging to $a$ and such that the
planes tangent to $\tau$ at $b_i$ (considered as the points in the associated
Grassmann bundle $G^{\dim \tau}(T_* \C^n)$) converge to a plane in $T_a \C^n$,
this limit plane does not lie in the hyperplane $\{df|_a = 0\} \subset
T_a\C^n$. If this condition is satisfied, we say that the one-parametric family
of varieties $A \cup f^{-1}(t),$ $t \in \C,$ has a {\it stratified Morse
singularity} at $a$.
\medskip

{\bf Notation.} For any $t \in \C$ denote by $X_t$ the set $f^{-1}(t)$ and by
$AX_t$ the set $ A\cup X_t$.
\medskip

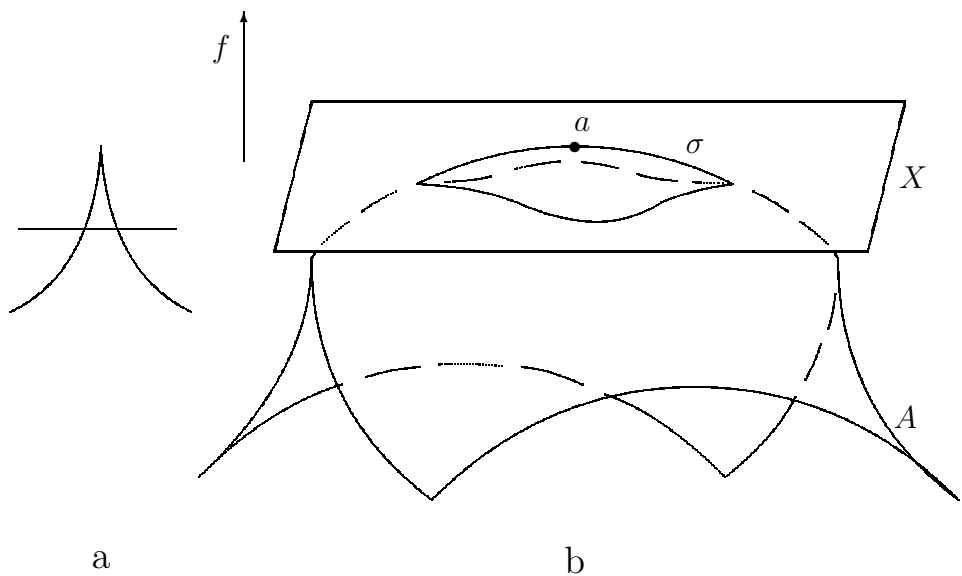
\begin{figure}
\unitlength 1.00mm \linethickness{0.4pt}
\begin{center}
\begin{picture}(126.00,74.00)
\bezier{180}(40.00,41.00)(40.00,21.00)(56.00,9.00)
\bezier{144}(40.00,41.00)(40.00,26.00)(25.00,12.00)
\bezier{180}(110.00,41.00)(110.00,21.00)(126.00,9.00)
\bezier{180}(56.00,9.00)(71.00,24.00)(91.00,24.00)
\bezier{180}(126.00,9.00)(111.00,24.00)(91.00,24.00)
\bezier{60}(95.00,12.00)(101.67,17.33)(104.33,22.67)
\bezier{20}(105.67,25.00)(107.00,27.33)(107.67,29.67)
\bezier{20}(108.67,31.67)(109.33,33.67)(109.67,36.33)
\bezier{80}(95.00,12.00)(86.33,20.67)(79.00,23.00)
\bezier{36}(76.00,24.33)(71.67,26.33)(68.00,26.67)
\bezier{20}(65.33,26.80)(61.00,27.20)(57.00,27.00)
\bezier{28}(54.33,26.67)(50.33,26.33)(47.33,25.33)
\bezier{92}(44.00,24.33)(34.33,21.33)(25.00,12.00)
\put(75.00,56.00){\circle*{1.33}} \put(76.00,59.00){\makebox(0,0)[cc]{$a$}}
\put(91.00,56.00){\makebox(0,0)[cc]{$\sigma$}}
\put(119.00,20.00){\makebox(0,0)[cc]{$A$}}
\put(35.00,42.00){\line(1,0){79.00}}
\multiput(114.00,42.00)(0.12,0.48){42}{\line(0,1){0.48}}
\put(119.00,62.00){\line(-1,0){79.00}}
\multiput(40.00,62.00)(-0.12,-0.48){42}{\line(0,-1){0.48}}
\bezier{36}(86.67,48.67)(82.67,46.00)(78.00,46.00)
\bezier{40}(86.67,48.67)(91.00,50.33)(96.00,51.00)
\bezier{48}(78.00,46.00)(73.00,46.00)(67.00,48.67)
\bezier{52}(67.00,48.67)(61.67,50.67)(54.00,51.00)
\bezier{88}(96.00,51.00)(86.00,56.00)(75.00,56.00)
\bezier{88}(75.00,56.00)(64.00,56.00)(54.00,51.00)
\bezier{16}(51.25,50.00)(48.33,48.67)(46.67,47.33)
\bezier{12}(45.00,46.00)(43.50,45.00)(41.33,43.00)
\bezier{16}(108.67,43.00)(106.50,45.00)(105.00,46.00)
\bezier{20}(103.33,47.33)(101.67,48.67)(98.75,50.00)
\multiput(40.00,41.00)(0.11,0.17){6}{\line(0,1){0.17}}
\multiput(110.00,41.00)(-0.11,0.17){6}{\line(0,1){0.17}}
\bezier{24}(63.80,52.00)(61.80,51.50)(57.80,51.25)
\bezier{20}(73.00,54.00)(70.50,53.80)(67.00,53.00)
\bezier{20}(76.00,53.90)(78.67,53.70)(80.67,53.10)
\bezier{24}(83.33,52.30)(85.33,51.70)(88.50,51.43)
\bezier{12}(90.67,51.25)(91.80,51.15)(96.00,51.05)
\bezier{116}(12.00,56.00)(11.00,39.33)(0.00,34.00)
\bezier{116}(12.00,56.00)(13.00,39.33)(24.00,34.00)
\put(1.00,45.00){\line(1,0){21.00}}
\put(12.00,1.00){\makebox(0,0)[cc]{{\large a}}}
\put(75.00,1.00){\makebox(0,0)[cc]{{\large b}}}
\put(31.00,54.00){\vector(0,1){20.00}}
\put(28.00,69.00){\makebox(0,0)[cc]{$f$}}
\put(120.00,52.00){\makebox(0,0)[cc]{$X$}}
\end{picture}
\caption{A Morse function on a stratified variety}
\end{center}
\end{figure}

For example, in Fig.~1b (the complexification of) the linear function defining
the plane $X$ has a stratified Morse singularity on (the complexification of)
$A$.
\medskip

Let $L_{\al}$ be the linear (i.e. with fibre $\C^1$) local system on $\C^n
\setminus AX_t$ with the set $\al = (\al_1, \ldots, \al_\nu)$ of ramification
indices: the cardinality $\nu$ of this set is equal to the number of
irreducible $(n-1)$-dimensional components of $AX_t$, and any small loop going
around a smooth piece of the $i$-th component in the positive direction (with
respect to the natural complex structure in the normal bundle) acts on the
fibre as multiplication by $\al_i$. We assume that $X_t$ is the first component
of $AX_t,$ so that $\al_1$ is responsible for the rotations around $X_t$.

The dual  local system will be denoted by $L_{\al^*},$ $\al^* = (\al_1^{-1},
\ldots, \al_\nu^{-1})$. For all $t \in \C \setminus 0$ sufficiently close to
$0$ (say, satisfying the condition $|t| \le \delta$ with sufficiently small
$\delta$) all corresponding sets $\C^n \setminus AX_t$ are homeomorphic to one
another and form a locally trivial bundle over the set of such $t$.

Denote by $C$ the loop $\{\delta e^{i\tau} \}, \ \tau \in [0,2\pi],$ generating
the group $\pi_1(\C^1 \setminus 0)$; the inverse loop $\{\delta e^{-i\tau}\}$
is denoted by $C^*$, and the closed disc bounded by any of these loops by
$D_{\delta}$.

The monodromy action of the loop $C$ (and $C^*$ as well) in two groups
$$ H_*(\C^n \setminus AX_t, \La) \quad and \quad
H^{lf}_*(\C^n \setminus AX_t, \La)$$ is well defined, where the letters $lf$ in
the last expression denote the homology groups of locally finite chains. The
study of this action is our main aim.
\medskip

These monodromy operators can be localized in the standard way, see e.g.
\cite{Ph2}, \cite{V1}. We suppose that $a$ is the unique point of
non-transversality of the manifold $X_0$ and the stratified variety $A$, and
the number $\delta$ participating in the definition of the loop $C$ is so small
that for any $t \in D_{\delta} \setminus 0$ the manifold $X_t$ is smooth and
transversal to the naturally stratified variety $A \cup \partial B $.

We choose $\{\delta\}$ as the distinguished point in $D_{\delta}\setminus 0$
and redenote $X_{\delta}$ and $AX_{\delta}$ simply by $X$ and $AX$. For any
non-negative integer $i $ we consider four groups
\begin{equation}
\label{groups}
\begin{array}{rclrcl}
\bar {\mathcal H}_{i,\al} & \equiv & H_i^{lf}(B \setminus AX,\partial B;L_\al),
\quad & {\mathcal H}_{i,\al} & \equiv & H_i^{lf}(B \setminus AX,L_\al), \\
\bar \chi_{i,\al} & \equiv & H_i(B \setminus AX,\partial B;L_\al), \quad &
\chi_{i,\al} & \equiv & H_i(B \setminus AX,L_\al).
\end{array}
\end{equation}

There are obvious homomorphisms
\begin{equation}
\label{triple}
\begin{array}{rclrcl}
\tilde i_{\al}: {\mathcal H}_{*,\al} & \to & H^{lf}_*(\C^n \setminus AX, \La),
& \quad \tilde j_{\al}: H^{lf}_*(\C^n \setminus AX, \La) & \to & \bar {\mathcal H}_{*,\al}, \\
i_{\al}: \chi_{*,\al} & \to & H_*(\C^n \setminus AX, \La), & \quad j_{\al}:
H_*(\C^n \setminus AX, \La) & \to & \bar \chi_{*,\al}.
\end{array}
\end{equation}

Also, any loop $\lambda \in \pi_1(D_\delta \setminus 0)$ defines in the
standard way (see \cite{AVGL}, \cite{V1}) the {\em local variation operators}
\begin{equation}
\label{varrs} \widetilde{Var}_{(\lambda)}: \bar {\mathcal H}_{*,\al} \to
{\mathcal H}_{*,\al} \, , \quad Var_{(\lambda)}: \bar \chi_{*,\al} \to
\chi_{*,\al} \ ,
\end{equation}
so that  the monodromy action of the loop $\lambda$ in the group $
H^{lf}_*(\C^n \setminus AX, \La) $ (respectively, $ H_*(\C^n \setminus AX,
\La)$) is equal to $\hbox{Id} + \tilde i_{\al}\circ \widetilde{Var}_{(\lambda)}
\circ \tilde j_{\al}$ (respectively, $\hbox{Id} + i_{\al}\circ Var_{(\lambda)}
\circ j_{\al}). $
\medskip

The composition operators $\tilde J_\alpha \equiv \tilde j_\alpha \circ \tilde
i_\alpha : {\mathcal H}_{*,\al} \to \bar {\mathcal H}_{*,\al}$ and $J_\alpha
\equiv j_\alpha \circ i_\alpha : \chi_{*,\al} \to \bar \chi_{*,\al}$ allow us
to express the {\it local monodromy} action of the loop $\lambda$ on four
groups $\bar {\mathcal H}_{*,\alpha}, \bar \chi_{*,\alpha}, {\mathcal
H}_{*,\alpha}, \chi_{*,\alpha}$ as
\begin{equation}
\label{mons}
\begin{array}{rr}
Id+ \tilde J_\alpha \circ \widetilde{Var}_{(\lambda)}, \quad &
Id+ J_\alpha \circ Var_{(\lambda)}, \\
Id+ \widetilde{Var}_{(\lambda)} \circ \tilde J_\alpha, \quad & Id+
Var_{(\lambda)} \circ J_\alpha
\end{array}
\end{equation}
respectively.

The groups (\ref{groups}) are related by non-degenerate Poincar\'e--Lefschetz
pairings
\begin{equation}
\label{poinc}
\begin{array}{rcl}
\bar {\mathcal H}_{i,\al} \otimes \chi_{2n-i,\al^*} & \to & \C ,  \\
{\mathcal H}_{i,\al} \otimes \bar \chi_{2n-i,\al^*} & \to & \C.
\end{array}
\end{equation}

{\bf Proposition 1.} {\it For any $i$ and $\alpha$, the operators
\begin{equation}
\label{vars}
\begin{array}{rccl}
\widetilde{Var}_{(C)}: & \bar {\mathcal H}_{i,\al} & \to & {\mathcal H}_{i,\al} \ , \\
Var_{(C^*)}: & \bar \chi_{2n-i,\al^*} & \to & \chi_{2n-i,\al^*}
\end{array}
\end{equation}
are conjugate with respect to the pairings (\ref{poinc}), i.e., for any
elements $ x \in \bar {\mathcal H}_{i,\al} $ and $ y \in \bar \chi_{2n-i,\al^*}
$ we have}
$$ \langle x, Var_{(C^*)}y \rangle =
\langle \widetilde{Var}_{(C)} x, y \rangle.$$

Proof (cf. \cite{sabir}). Let $x$ and $y$ be some elements of the groups $ \bar
{\mathcal H}_{i,\al} $ and $  \bar \chi_{2n-i,\al^*} $ respectively, and ${\bf
x}$ and ${\bf y}$ some relative cycles representing them and intersecting
generically in $B$ (in particular having no common points in $\partial B$).
Then the intersection number $\langle{\bf x},{\bf y}\rangle$ is well defined
(although it is not an invariant of the homology classes $x$ and $y$). Let
$h_C{\bf x}$ and $h_{C^*}{\bf y}$ be two relative cycles in $B \setminus AX$
obtained from ${\bf x}$ and ${\bf y}$ by the monodromy over the loops $C$ and
$C^*$ as in the construction of variation operators (i.e., fixed close to the
boundary of $B$). Then $\langle\widetilde{Var}_{(C)}(x),y\rangle - \langle x,
Var_{(C^*)}(y)\rangle = \langle h_C {\bf x} - {\bf x},{\bf y}\rangle - \langle
{\bf x},h_{C^*} {\bf y} - {\bf y} \rangle = \langle h_C{\bf x}, {\bf y} \rangle
- \langle {\bf x}, h_{C^*} {\bf y} \rangle = \langle h_C{\bf x}, {\bf y}
\rangle - \langle h_{C^*}h_C{\bf x}, h_{C^*} {\bf y} \rangle =0$ because these
intersection pairings are invariant under the monodromy action. \quad $\Box$
\medskip

Set $k=\dim \sigma$, $m=n-k,$ and let $T$ be the $m$-dimensional complex plane
through $a$ transversal to $\sigma$. Then, replacing the sets $B, A$ and $X_t$
by their intersections with $T,$ we get all structures as above (homology
groups, local variation and monodromy action, intersection indices and maps
$\tilde J_\al,$ $J_\al$) placed in the space of reduced dimension $m$. In the
next section we show how, knowing all these objects for this reduced space, we
can reconstruct them for the initial objects in entire space $\C^n.$ This
reconstruction will de done by induction over a flag of planes connecting $T$
and $\C^n.$
\medskip

\section{Main results}

Consider the flag of complex planes (or, more generally, smooth complex
submanifolds)
\begin{equation} \label{flag}
T \equiv T^m \subset T^{m+1} \subset \cdots \subset T^n \equiv \C^n,
\end{equation}
of dimensions $m, m+1, \ldots, n$ respectively, all of which intersect $\sigma$
transversally at the point $a$. Then the intersections with the strata of $A$
define a Whitney stratification on any of varieties $A \cap B \cap T^r$. If the
flag (\ref{flag}) is generic, then for any $r = m, \ldots, n$ the restriction
of $f$ on $A \cap T^r$ is a Morse function in $B$; we shall suppose that this
condition is satisfied for all $r$. Any loop $\lambda$ in $D_\delta \setminus
0$, considered as family of planes $X_t \cap T^r$, $t \in \lambda$, defines
then the variation operators
\begin{equation}\label{varr}
\widetilde{\hbox{Var}}_{(\lambda),r}: H_i^{lf}(B \cap T^r \setminus AX,\partial
B;L_\al) \to H_i^{lf}(B \cap T^r \setminus AX,L_\al),
\end{equation}
\begin{equation}\label{varrr}
\hbox{Var}_{(\lambda),r}: H_i(B \cap T^r \setminus AX,\partial B;L_\al) \to
H_i(B \cap T^r \setminus AX,L_\al),
\end{equation}
in particular $\widetilde{\hbox{Var}}_{(\lambda),n} =
\widetilde{\hbox{Var}}_{(\lambda)}$, $\hbox{Var}_{(\lambda),n} =
\hbox{Var}_{(\lambda)}$.

Everywhere below we consider only operators (\ref{varr}), (\ref{varrr}) defined
by the loop $\lambda = C$ and denote these operators simply by $\hbox{Var}_r$.

Denote the groups participating in these operators as follows:
\begin{equation}
\label{grflag}
\begin{array}{ll}
\bar {\mathcal H}_{i,\al}(r)  \equiv H_i^{lf}(B \cap T^r \setminus AX,\partial
B;L_\al), &
{\mathcal H}_{i,\al}(r)  \equiv  H_i^{lf}(B \cap T^r \setminus AX,L_\al), \\
\bar \chi_{i,\al}(r)  \equiv H_i(B \cap T^r \setminus AX,\partial B;L_\al), &
\chi_{i,\al}(r) \equiv  H_i(B \cap T^r \setminus AX,L_\al),
\end{array}
\end{equation}
and local monodromy operators defined by the loop $C$ on four groups $\bar
{\mathcal H}_{*,\alpha}(r),$ ${\mathcal H}_{*,\alpha}(r), $ $\bar
\chi_{*,\alpha}(r), $ $\chi_{*,\alpha}(r)$ by $\bar M_r, M_r, \bar \mu_r,
\mu_r$ respectively.

Of course, the analogues of the Poincar\'e duality (\ref{poinc}) and
Proposition 1 are valid for any $r$.
\bigskip

{\bf Theorem 1.} {\it For any $l=0, \ldots, n-m,$ and any $\al,$ there are
almost canonical isomorphisms}
\begin{equation}
\label{stab}
\begin{array}{rccl}
\Sigma^l : &{\mathcal H}_{i,\al}(m) & \to & {\mathcal H}_{i+l,\al}(m+l), \\
\dar^l : & \bar \chi_{i,\al}(m) & \to & \bar \chi_{i+l,\al}(m+l).
\end{array}
\end{equation}

The construction of these maps will be described in the next section. For the
explanation of the word ``almost'' see Remark in \S \ 3.3.
\medskip

The stabilization of groups $\bar {\mathcal H}_{i+l,\al}(m+l)$ and
$\chi_{i+l,\al}(m+l)$ depends very much on the fact, is the coefficient $\al_1$
(responsible for the ramification of $L_{\al}$ close to the plane $X_t$) equal
to 1 or not. To unify the corresponding formulae, set $\Phi_{i,\al}\equiv
\psi_{i,\al} \equiv 0$ if $\al_1 \ne 1;$ if $\al_1=1,$ we consider the local
system $L_{\hat \al}$ on $B \setminus A,$ isomorphic to $L_{\al}$ in a small
closed disc $ b \subset B \setminus X_t$ with center $a$ (i.e. having the same
monodromy coefficients $\al_2, \ldots, \al_\nu$ related to all components of
$A$) and set $\Phi_{i,\al}\equiv H_i^{lf}(B \cap T^m\setminus A, \partial B;
L_{\hat \al}) \simeq H_{i+2l}^{lf}(B \cap T^{m+l}\setminus A, \partial B;
L_{\hat \al}),$ $\psi_{i,\al}\equiv H_i(B \cap T^m\setminus A, L_{\hat \al})
\simeq H_i(B \cap T^{m+l}\setminus A, L_{\hat \al}).$
\bigskip

{\bf Theorem 1$'$.} {\it For any $l=1, \ldots, n-m,$ and any $\al,$ there are
almost canonical isomorphisms}
\begin{equation}
\label{stab2}
\begin{array}{rcl}
\bar {\mathcal H}_{i+l,\al}(m+l) & \simeq &  {\mathcal H}_{i,\al}(m)
\oplus \Phi_{i-l,\al} \oplus \Phi_{i-l+1,\al} , \\
\chi_{i+l,\al}(m+l) & \simeq & \bar \chi_{i,\al}(m)
\oplus \psi_{i+l,\al} \oplus \psi_{i+l-1,\al}. \\
\end{array}
\end{equation}

For any  $l=1, \ldots, n-m,$ we denote by $\rlh_l$ (respectively, $\infty_l$)
the injection ${\mathcal H}_{i,\al}(m) \to \bar {\mathcal H}_{i+l,\al}(m+l)$
(respectively, $\bar \chi_{i,\al}(m) \to \chi_{i+l,\al}(m+l)$) defined by the
first summands in the right-hand parts of (\ref{stab2}).
\bigskip

{\bf Theorem 2.} {\it For any $l = 1, \ldots, n-m,$

a) the obvious homomorphism $ \tilde J_{\al,m+l} \equiv \tilde j_{\al}\circ
\tilde i_{\al}: {\mathcal H}_{i+l,\al}(m+l) \to \bar {\mathcal
H}_{i+l,\al}(m+l)$ (i.e., the reduction mod $\partial B$) maps any element
$\Sigma^l(x),$ $x \in {\mathcal H}_{i,\al}(m),$ to
\smallskip

$\rlh_{l} \circ \widetilde{Var}_m \circ \tilde J_{\al,m}(x) \equiv \, \rlh_{l}
\circ (M_m - Id)(x)$ \quad if $l$ is even \quad and to
\smallskip

$- \rlh_{l} (2x + \widetilde{Var}_m \circ \tilde J_{\al,m}(x)) \equiv -
\rlh_{l} \circ (M_m + Id)(x)$ \quad if $l$ is odd.
\smallskip

b) the similar homomorphism $ J_{\al,m+l} \equiv j_{\al}\circ i_{\al}:
\chi_{i+l,\al}(m+l) \to \bar \chi_{i+l,\al}(m+l)$ is equal to zero on the last
two summands in (\ref{stab2}) and maps any element $\infty_l(y)$, $y \in \bar
\chi_{i,\al}(m),$ to
\smallskip

$\dar^l \circ J_{\al,m} \circ Var_m(y) \equiv \, \dar^l \circ (\bar
\mu_m-Id)(y) $ \quad if $l$ is even \quad and to
\smallskip

$-\dar^l (2y + J_{\al,m}\circ Var_m(y)) \equiv - \dar^l \circ (\bar \mu_m +
Id)(y) $ \quad if $l$ is odd.}
\bigskip

{\bf Theorem 3.} {\it For any $\al$ and any $l = 1, \ldots, n-m,$

a) the operator $\widetilde{Var}_{m+l}: \bar {\mathcal H}_{i+l,\al}(m+l) \to
{\mathcal H}_{i+l,\al}(m+l)$ is equal to zero on two last summands in the first
row of (\ref{stab2}) and maps the element $\rlh_l(x)$, $x \in {\mathcal
H}_{i,\al}(m),$ of the first summand to $\Sigma^l(x)$;

b) the operator $Var_{m+l}: \bar \chi_{i+l,\al}(m+l) \to \chi_{i+l,\al}(m+l)$
maps any element $\dar^l(y),$ $y \in \bar \chi_{i,\al}(m),$ to $\infty_l(y)$.}
\bigskip

{\bf Corollary.} {\it The local monodromy operators $M_{m+l}, \bar M_{m+l},
\bar \mu_{m+l}$ and $\mu_{m+l}$ respectively map the elements $\Sigma^l(x),
\rlh_l(x), \dar^l(y)$ and $\infty_l(y)$ (where $x \in {\mathcal
H}_{*,\alpha}(m),$ $ y \in \bar \chi_{*,\alpha}(m)$) into the elements $(-1)^l
\Sigma^l(M_m(x))$, $(-1)^l \rlh_l(M_m(x))$, $(-1)^l \dar^l(\bar \mu_m(y))$ and
$(-1)^l \infty_l(\bar \mu_m(y))$ respectively. Operators $\bar M_{m+l}$ and
$\mu_{m+l}$ act trivially (i.e. as the identity operators) on the last two
summands in both formulae (\ref{stab2}).}
\medskip

{\bf Theorem 4.} {\it For any $l=1, \ldots, n-m,$

\noindent A) for any $x \in {\mathcal H}_{i,\al}(m)$ and $y^* \in \bar
\chi_{i,\al^*}(m),$
$$\langle \Sigma^l(x), \dar^l(y^*) \rangle =
(-1)^{l(i+1+(l-1)/2)}\langle x,y^* \rangle ; $$

\noindent B) for any $\xi \in \bar {\mathcal H}_{i+l,\al}(m+l) $ and $ \zeta^*
\in \chi_{2m-i+l,\al^*}(m+l)$,
\begin{enumerate}
\item
$\langle \xi,\zeta^* \rangle =0$ if $\xi$ and $\zeta^*$ belong to summands in
the right-hand part of (\ref{stab2}) placed not one over the other;
\item If $\xi = \rlh_l(x),$
$x \in {\mathcal H}_{i,\al}(m)$, and $\zeta^* = \infty_l(y^*)$, $ y^* \in \bar
\chi_{i,\al^*}(m)$, then $\langle \xi,\zeta^* \rangle= (-1)^{il+1+l(l-1)/2}
\langle x, \bar \mu_m(y^*) \rangle .$
\end{enumerate}
}
\medskip

{\bf Corollary} (periodicity theorem). {\it Suppose that the groups
$\Phi_{*,\al},$ $\psi_{*,\al}$ are trivial (e.g. $\al_1 \ne 1$). Then all maps
$\rlh_l,$ $\infty_l$ are isomorphisms, and the entire structure depending on
$A,$ $\al$ and $r$ and consisting of four graded groups $\bar {\mathcal
H}_{*,\al}(r),$ ${\mathcal H}_{*,\al}(r),$ $\bar \chi_{*,\al}(r),$
$\chi_{*,\al}(r),$ similar four groups with $\al$ replaced by $\al^*,$ all
possible intersection pairings between them, and the operators $\tilde J_\al,$
$J_\al,$ $\widetilde{Var}_r,$ $Var_r$, $M_{r}, \bar M_{r}, \mu_{r}$ and $\bar
\mu_{r}$ is periodic in $r$ with period 2.}
\bigskip

The ideas of proofs of some of these results are essentially contained in
\cite{V2}, however the {\it answers} are sometimes different and should have
been written explicitly: this it the main purpose of the present article.
\medskip

We shall often use the following well-known fact (see e.g. \cite{GM},
\cite{M}).
\medskip

{\bf Proposition 2.} {\it For any point $a$ of an analytic subset $A \subset
\C^n$ there exists a small disc $B$ centred at $a$ such that the pair $(A \cap
B,  A \cap \partial B)$ is homeomorphic to the cone over $A \cap \partial B$,
and this homeomorphism is identical on $A \cap \partial B$ and maps $a$ into
the vertex of the cone. Moreover, the same is true for any smaller concentric
disc.} $\quad \Box$

\section{Realization of formulae (\protect\ref{stab}), (\protect\ref{stab2})}

The stabilization of groups (\ref{grflag}) will be constructed by induction
over the flag of planes $T^{m+l}$. Namely, for any $r = m, \ldots ,n-1$ we
construct the maps
\begin{equation}
\label{stab3}
\begin{array}{rccl}
\Sigma : & {\mathcal H}_{j,\al}(r) & \to & {\mathcal H}_{j+1,\al}(r+1), \\
\dar : & \bar \chi_{j,\al}(r) & \to & \bar \chi_{j+1,\al}(r+1), \\
\rlh : & {\mathcal H}_{j,\al}(r) & \to & \bar {\mathcal H}_{j+1,\al}(r+1), \\
\infty : & \bar \chi_{j,\al}(r) & \to & \chi_{j+1,\al}(r+1).
\end{array}
\end{equation}
Two first of them (and also two last if $\Phi_{*,\alpha} = \psi_{*,\alpha}=0$)
are isomorphisms.
\medskip

Here are some preliminary constructions and reductions, cf. \cite{V2}.

\subsection{Adapted coordinates and polydisk $B'$}

\noindent Let $A, a, \sigma, B, f, D_\delta$ and $X_t$ be the same as in the
previous sections. \smallskip

{\bf Definition.} A local analytic coordinate system $\{z_1, \ldots, z_n\}$ in
$\C^n$ with origin at $a$ is called {\it adapted} if $z_n \equiv f$, the
tangent space to $\sigma$ at $a$ is spanned by the vectors $\partial/\partial
z_1, \ldots, \partial/\partial z_k$ (so that the restrictions of the functions
$z_1, \ldots, z_k$ constitute a local coordinate system on $\sigma$), and in
restriction to $\sigma$
$$z_n \equiv z_1^2 + \cdots + z_k^2 . $$

In Fig. 1b a real version of this situation is shown, where $n=3$, the plane $X
\equiv X_\delta$ is given by the equation $z_3 \equiv \delta$, and $z_1$ is the
coordinate along the stratum $\sigma$. The transversal slice of this picture by
the plane $\{z_1 = 0\}$ is shown in Fig. 1a. \smallskip

By the Morse lemma, adapted coordinates always exist; let us fix such a
coordinate system. Without loss of generality we can define the flag
(\ref{flag}) by the conditions $T^r = \{z \mid \nolinebreak z_1 = \cdots =
z_{n-r} =0\}$.

We can assume that $B \equiv B_\varepsilon$ is a closed disc of radius
$\varepsilon$ with respect to the standard Hermitian metric defined by these
coordinates $z_1, \ldots, z_n$. Moreover, in our considerations we can replace
it by a closed polydisk defined by these coordinates. Namely, let $B' \subset
B$ be the polydisk $\{z \mid |z_i| \le \varepsilon/n \hbox{ for all } i\}$ and
suppose that the number $\delta$ (participating in the definition of the disc
$D \equiv D_\delta \subset \C^1$) is sufficiently small with respect to
$\varepsilon$ and $\varepsilon^2$.
\medskip

{\bf Lemma 1.} {\it For every surface $X_\lambda \subset B$ defined by the
equation $z_n \equiv \lambda, \ \lambda \in D$, the pair $(B' \setminus
AX_\lambda, \partial B' \setminus AX_\lambda)$ is homeomorphic to the pair $(B
\setminus AX_\lambda, \partial B \setminus AX_\lambda).$ These homeomorphisms
depend continuously on the parameter $\lambda$, and the induced isomorphisms of
all groups (\ref{groups}) onto similar groups in whose definition $B$ is
replaced by $B'$ coincide with the morphisms induced by the identical embedding
$B' \setminus AX_\lambda \to B \setminus AX_\lambda$ in the case of groups
${\mathcal H}_{i,\al}, \chi_{i,\al}$ and to the morphisms induced by the
obvious map $(B \setminus AX_\lambda)/(\partial B \setminus AX_\lambda) \to (B'
\setminus AX_\lambda)/(\partial B' \setminus AX_\lambda)$ (reduction mod $B
\setminus B'$) in the case of groups $\bar {\mathcal H}_{i,\al}, \bar
\chi_{i,\al}.$

Moreover, for every $r = m, \ldots, n$ and every $\lambda \in D$, all these
statements remain valid if we replace both $B$ and $B'$ by $B \cap T^r$ and $B'
\cap T^r$ respectively.}
\bigskip

The proof repeats the proof of Proposition 5 in \cite{V2}.
\medskip

Thus, everywhere in the proof of Theorems 1--4 we can replace the disc $B$ by
the polydisk $B'$.

\subsection{Fibre bundle \protect$z_{n-r}$}

\noindent For any $r = m, \ldots, n$ denote the $2r$-dimensional polydisk $B'
\cap T^r$ by $B^{(r)}$ and the characteristic radius $\varepsilon /n$ of all
such polydisks by $\epsilon$. Remember the notation $X \equiv X_\delta$ and $AX
\equiv A \cup X_\delta$. Denote the $\epsilon$-disc $\{w \in \C^1 \mid |w| \le
\epsilon \}$ by $\Omega$. For arbitrary $r < n$ consider the projection
\begin{equation} \label{przet}
z_{n-r}: (B^{(r+1)} \setminus AX)  \to \Omega.
\end{equation}
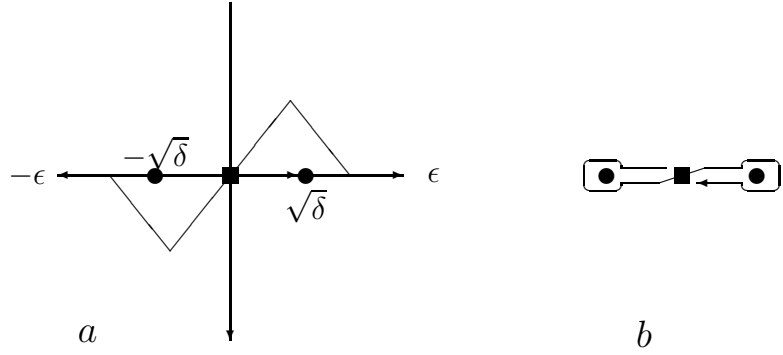
\begin{figure}
\unitlength=1.00mm \special{em:linewidth 0.4pt} \linethickness{0.4pt}
\begin{picture}(113.00,53.00)
\put(39.00,29.00){\rule{2.00\unitlength}{2.00\unitlength}}
\put(30.00,30.00){\circle*{2.00}} \put(50.00,30.00){\circle*{2.00}}
\put(51.00,30.00){\vector(1,0){12.00}} \put(29.00,30.00){\vector(-1,0){12.00}}
\put(31.00,30.00){\vector(1,0){18.00}} \put(40.00,53.00){\vector(0,-1){45.00}}
\put(67.00,30.00){\makebox(0,0)[cc]{$\epsilon$}}
\put(13.00,30.00){\makebox(0,0)[cc]{$-\epsilon$}}
\put(30.00,33.00){\makebox(0,0)[cc]{$-\sqrt{\delta}$}}
\put(50.00,26.00){\makebox(0,0)[cc]{$\sqrt{\delta}$}}
\put(24.00,30.00){\line(4,-5){8.00}} \put(32.00,20.00){\line(4,5){16.00}}
\put(48.00,40.00){\line(4,-5){8.00}} \put(95.00,9.00){\makebox(0,0)[cc]{{\large
$b$}}} \put(21.00,9.00){\makebox(0,0)[cc]{{\large $a$}}}
\put(99.00,29.00){\rule{2.00\unitlength}{2.00\unitlength}}
\put(110.00,30.00){\circle*{2.00}} \put(90.00,30.00){\circle*{2.00}}
\put(98.00,31.00){\line(-1,0){6.00}} \put(92.00,29.00){\line(1,0){5.00}}
\put(97.00,29.00){\line(3,1){6.00}} \put(103.00,31.00){\line(1,0){5.00}}
\put(108.00,29.00){\vector(-1,0){6.00}} \put(110.50,31.00){\oval(5.00,2.00)[t]}
\put(110.50,29.00){\oval(5.00,2.00)[b]} \put(89.50,29.00){\oval(5.00,2.00)[b]}
\put(89.50,31.00){\oval(5.00,2.00)[t]} \put(87.00,31.00){\line(0,-1){2.00}}
\put(113.00,29.00){\line(0,1){2.00}}
\end{picture}
\caption{Base \protect$\Omega$ of the fibre bundle
\protect$z_{n-r}$}\label{omega}
\end{figure}

For any $t \in \Omega$ denote by ${\mathcal F}_t$ the fibre $(B^{(r+1)}
\setminus AX ) \cap \{z|z_{n-r}=t\}$ of this projection over the point $t$, and
by $\partial {\mathcal F}_t$ the boundary ${\mathcal F}_t \cap
\{z|\max(|z_{n-r+1}|, \ldots, |z_n|)=\delta\} \subset \partial B^{(r+1)}$ of
this fibre. In particular, ${\mathcal F}_0 \equiv B^{(r)}\setminus AX$,
$H_*^{lf}({\mathcal F}_0,
\partial {\mathcal F}_0; L_\alpha) \simeq \tilde {\mathcal H}_{*, \alpha}(r),$
$H_*^{lf}({\mathcal F}_0, L_\alpha) \simeq {\mathcal H}_{*, \alpha}(r),$
$H_*({\mathcal F}_0, \partial {\mathcal F}_0; L_\alpha) \simeq \tilde \chi_{*,
\alpha}(r),$ and $H_*({\mathcal F}_0, L_\alpha) \simeq \chi_{*, \alpha}(r).$

Any of two planes $\{z|z_{n-r}= \pm \sqrt{\delta}\}$ contains a distinguished
point $\Delta^{\pm}$, the critical point of the restriction of $z_{n-r}$ to the
manifold $\sigma \cap X \cap T^{r+1}$.
\smallskip

{\bf Lemma 2.} {\it If the radius $\epsilon$ of the polydisk $B'$ is
sufficiently small and $\delta$ is sufficiently small with respect to
$\epsilon^2$, then

a) the projection (\ref{przet}) defines a locally trivial fibre bundle over the
disc $\Omega$ with two points $\pm \sqrt{\delta}$ removed, the standard fibre
$({\mathcal F}_t, \partial {\mathcal F}_t)$ of which is homeomorphic to
$({\mathcal F}_0,
\partial {\mathcal F}_0)$;

b) the fibres over the exceptional points $\sqrt{\delta}$ and $-\sqrt{\delta}$
are homeomorphic to the direct product $\partial {\mathcal F}_0 \times (0,1]$
(or, more transparently, to the disc $B^{(r)}$ from which a cone over $\partial
B^{(r)} \setminus \partial {\mathcal F}_0$ is removed);

c) the restriction of this projection on $\partial B^{(r+1)}$ defines a
trivializable bundle with typical fibre $\partial B^{(r)} \setminus AX$ over
the interior part of $\Omega$.}
\smallskip

This lemma follows directly from the construction, from Thom's isotopy theorem
and from Proposition 2, cf. \cite{V2}. $\quad \Box$
\smallskip

Consider two variation operators $\tilde V_{+,-}: H^{lf}_*({\mathcal F}_0,
\partial {\mathcal F}_0;L_\al) \to H^{lf}_*({\mathcal F}_0,L_\al)$ and two similar
operators $V_{+,-}: H_*({\mathcal F}_0, \partial {\mathcal F}_0;L_\al) \to
H_*({\mathcal F}_0,L_\al)$ defined by the simple loops in $\Omega$
corresponding to the segments $[0, \sqrt{\delta}]$ and $[0, -\sqrt{\delta}]$.
(For instance the $\infty$-shaped loop in Fig. \ref{omega}b is a composition of
the second of these loops and the loop inverse to the first of them.)
\smallskip

{\bf Lemma 3.} {\it The operators $\tilde V_+$, $\tilde V_-$ are equal to one
another and to the operator $\widetilde{Var}_r: \bar {\mathcal H}_*(r) \to
{\mathcal H}_*(r).$ The operators $V_+$, $V_-$ are equal to one another and to
the operator} $Var_r: \bar \chi_*(r) \to \chi_*(r).$
\medskip

Indeed, any of these operators is defined by a loop in the space of pairs of
complex hyperplanes in $T^{r+1}$ in general position with $A$: for $V_+$ and
$V_-$ the first plane of this pair is fixed and coincides with $X_\delta \cap
T^{r+1}$, while the second is distinguished by the condition $z_{n-r}=t$ where
$t$ runs over the corresponding simple loop in $\Omega$; for $Var_r$ the {\it
second} plane is fixed and coincides with the one distinguished by $z_{n-r}=0$,
and the first moves and coincides with the planes $X_\lambda \cap T^{r+1}$
where $\lambda$ runs over the circle $C \equiv \partial D_\delta$.

It is easy to see that all these three loops are homotopic to one another in
the space of pairs of hyperplanes generic with respect to $A$, and Lemma 3 is
proved. $\quad \Box$

\subsection{Construction of maps (\protect\ref{stab3}).}

Let $x$ be any element of the group ${\mathcal H}_{j,\al}(r)$. The class
$\Sigma(x) \in {\mathcal H}_{j+1,\al}(r+1)$ realizing the first isomorphism in
(\ref{stab3}) is obtained from $x$ by a sort of the suspension operation.
Namely, using the fibre bundle structure (\ref{przet}), we transport the
realizing $x$ cycle ${\bf x} \subset {\mathcal F}_0$ over the segment
$[-\sqrt{\delta}, \sqrt{\delta}]$. This one-parametric family of cycles sweeps
out a $(j+1)$-dimensional chain in $B^{(r+1)} \setminus AX$ oriented by the
pair of orientations, the first of which is the orientation of the base segment
chosen as in Fig. \ref{omega}a, and the second is induced by the original
orientation of ${\bf x}$. The boundary of this chain lies in the marginal
fibres over the endpoints $-\sqrt{\delta}$ and $\sqrt{\delta}$. By the
statement b) of Lemma 2 the homology groups with closed supports of these
fibres are trivial, thus we can contract these boundaries inside these fibres
and get a cycle; $\Sigma(x)$ is defined as its homology class.

To construct the cycle $\rlh(x)$ we first transport the cycle ${\bf x}$ over
the S-shaped path in Fig. \ref{omega}a and get similar cycles in the fibres
over some two interior points of the segments $[\sqrt{\delta}, \epsilon]$ and
$[-\epsilon, -\sqrt{\delta}]$, then sweep out some $(j+1)$-dimensional chains
over these segments oriented as is shown in the same Fig. \ref{omega}a, and
again span the boundaries of the obtained chains inside the fibres over the
points $\pm \sqrt{\delta}$. The class $\rightleftharpoons(x) \in \bar {\mathcal
H}_{j+1,\al}(r+1)$ corresponding to $x$ is defined by half the difference of
these two cycles.

Given a homology class $y \in \bar \chi_{j,\al}(r)$, the corresponding cycle
$\dar(y) \in \bar \chi_{j+1,\al}(r+1)$ is swept out by the similar
one-parametric family of cycles obtained from a realizing $y$ compact cycle
${\bf y}$ transported over the axis $\{Re\ z_{n-r} =0\}$ oriented downwards.

To obtain the class $\infty(y)  \in \chi_{j+1, \al}(r+1)$ we transport the same
cycle ${\bf y}$ along the $\infty$-shaped path in Fig. \ref{omega}b. The
boundary of the $(j+1)$-dimensional chain swept out by it belongs to $\partial
B'$ and is homeomorphic there (via the trivializing homeomorphism from
statement c) of Lemma 2) to the direct product of $\partial {\bf y}$ and this
path. In particular, it is the boundary of the $(j+1)$-dimensional chain in
$\partial B'$, homeomorphic to the direct product of $\partial {\bf y}$ and the
2-chain in $\C^1$ bounded by this path.  Thus the difference of these two
$(j+1)$-dimensional chains is an absolute compact cycle in $B^{(r+1)} \setminus
AX$; $\infty(y)$  is defined as its homology class.
\medskip

Operations $\Sigma^l$ and $\dar^l$ participating in Theorem 1 are just the
$l$-fold iterations of these homomorphisms $\Sigma$ and $\dar$. The
homomorphisms $\rlh_l$ and $\infty_l$ defining the first summands in
(\ref{stab2}) are defined as $\rlh \circ \Sigma^{l-1}$ and $\infty \circ
\dar^{l-1}$ respectively.
\medskip

{\bf Remark.} If we replace the coordinate $z_{n-r}$ by $-z_{n-r}$, all
homomorphisms $\Sigma,$ $\dar,$ $\rlh$ and $\infty$ will be multiplied by $-1$.
The operations $\Sigma^l,$ $\dar^l,$ $\rlh_l$ and $\infty_l$ are thus defined
only up to a (common) sign switching with any involution $z_q \to z_{-q},$
$q=k, k-1, \ldots, k-l+1$.
\medskip

{\bf Conjecture.} {\it There is a continuous involution of the pair
$(B^{(r+1)}, AX)$ commuting with the involution $z_{n-r} \to -z_{n-r}$ of
$\Omega$ and identical on the fibre $\{z|z_{n-r}=0\}$ over the fixed point of
the latter involution.}
\medskip

Although we do not prove this conjecture, in the next section we use its
``homological shadow'', i.e. the involution in the homology groups, existing
(as we shall see) independently on this conjecture.  \medskip

Finally, we realize two last summands in both equations (\ref{stab2}). Let
$\Delta \in B^{(m+l)}$ be any intersection point of the plane $X_{\delta}$ and
the stratum $\sigma$ of $A$. Let $\beta \subset B^{(m+l)}$ be a very small
closed disc centred at $\Delta$. Then by the K\"unneth formula for arbitrary
$\al$
\begin{equation}
\label{kunn} H_*(\beta \setminus AX, L_{\al}) \simeq \psi_{*,\al} \otimes
H_*(\C^1 \setminus 0, L_{\al_1}),
\end{equation}
\begin{equation}
\label{kunn1} H^{lf}_*(\beta \setminus AX, \partial \beta ; L_{\al}) \simeq
\Phi_{*-2l,\al} \otimes H^{lf}_*(\C^1 \setminus 0, L_{\al_1})
\end{equation}
(where the groups $H_*(\C^1 \setminus 0, L_{\al_1})$  and $H^{lf}_*(\C^1
\setminus 0, L_{\al_1})$ are trivial if $\al_1 \ne 1$ and are isomorphic to
$H_*(S^1)$ and $H_{*-1}(S^1)$ respectively if $\al_1 = 1$, and $\Phi_{*-2l}$
is the graded group obtained from $\Phi_*$ by the shift of grading). The sum of
two last summands in the second row of (\ref{stab2}) is isomorphic to the
$(i+l)$-dimensional component of the group (\ref{kunn}) and is realized as
follows.  \medskip

{\bf Lemma 4.} {\it The homomorphism
\begin{equation}
\label{inj} H_*(\beta \setminus AX, L_{\al}) \to \chi_{*,\al}(m+l)
\end{equation}
induced by the identical embedding is injective.}
\medskip

If $l>1$ (so that the set $B^{(m+l)} \cap X_\delta \cap \sigma$ of possible
points $\Delta$ is path-connected) then the sum of two last summands in the
second row of (\ref{stab2}) coincides with the image of this injection. If
$l=1$ and this set consists of two points $\Delta^+$ and $\Delta^-$ (see \S \
3.2), then the corresponding groups $H_*(\beta^+ \setminus AX, L_{\al})$,
$H_*(\beta^- \setminus AX, L_{\al})$ are naturally identified to one another.
(This identification is transparent on the dual cohomological level, because
both cohomology groups dual to two factors in (\ref{kunn}) are induced from
groups $H^*(B^{(m+l)} \setminus A)$ and $H^*(B^{(m+l)} \setminus X)$
respectively.)

Consider the similar homomorphism $H_*(\beta^+ \setminus AX, L_{\al}) \oplus
H_*(\beta^- \setminus AX, L_{\al}) \to \chi_{*,\al}(m+l)$ defined by the
embedding $(\beta^+ \cup \beta^-) \to B^{(m+l)}$ and take the subgroup of its
image invariant under this identification; this subgroup is again isomorphic to
either of $H_*(\beta^\pm \setminus AX, L_{\al})$  and realizes the last
summands in the second row of (\ref{stab2}). Two last summands in the first row
are Poincar\'e dual to them, let us describe them explicitly.

The pair $(B^{(m+l)}, A)$ is homeomorphic to the product $(B^{(m)}, B^{(m)}
\cap A) \times \C^l,$ therefore for any cycle $w \in \Phi_{i-l}$ the cycle $w
\times \C^l$ defines an element of the group $H^{lf}_{i+l}(B^{(m+l)} \setminus
A, \partial B; L_{\hat \alpha})$ and, since $\alpha_1 =1,$ also of the group
$H^{lf}_{i+l}(B^{(m+l)} \setminus AX, \partial B; L_{\alpha}) \equiv \bar
{\mathcal H}_{i+l,\al}(m+l) $. We choose this cycle in general position with a
cycle generating the group $H^{lf}_{2(m+l)-1}(B^{(m+l)} \setminus X, \partial
B)$ (say, with the cycle given by the condition $z_n < \delta$), then the
intersection of these two cycles gives us also an element of the group $\bar
{\mathcal H}_{i+l-1,\al}(m+l) $; this element corresponds to the element $z$ of
the third summand of the decomposition (\ref{stab2}) of this group.

\section{Proof of main theorems}

\begin{figure}
\unitlength=1.00mm \special{em:linewidth 0.4pt} \linethickness{0.4pt}
\begin{picture}(137.00,116.00)
\put(15.00,105.00){\circle*{2.00}} \put(25.00,105.00){\circle*{2.00}}
\put(20.00,105.00){\oval(30.00,20.00)[]} \put(39.00,105.00){\vector(1,0){7.00}}
\put(65.00,105.00){\oval(30.00,20.00)[]} \put(84.00,105.00){\vector(1,0){7.00}}
\put(110.00,105.00){\oval(30.00,20.00)[]} \put(115.00,45.00){\circle*{2.00}}
\put(105.00,45.00){\circle*{2.00}} \put(114.83,45.67){\oval(4.33,4.00)[t]}
\put(112.67,45.67){\line(-1,0){6.67}} \put(104.67,44.33){\oval(4.67,4.00)[b]}
\put(107.00,44.33){\line(1,0){7.00}} \put(20.00,75.00){\oval(30.00,20.00)[]}
\put(39.00,75.00){\vector(1,0){7.00}} \put(65.00,75.00){\oval(30.00,20.00)[]}
\put(84.00,75.00){\vector(1,0){7.00}} \put(110.00,75.00){\oval(30.00,20.00)[]}
\put(15.00,75.00){\circle*{2.00}} \put(25.00,75.00){\circle*{2.00}}
\put(20.00,75.00){\vector(-1,1){5.00}} \put(15.00,80.00){\vector(-1,-1){5.00}}
\put(10.00,75.00){\vector(1,-1){5.00}} \put(15.00,70.00){\vector(1,1){10.00}}
\put(25.00,80.00){\vector(1,-1){5.00}} \put(30.00,75.00){\vector(-1,-1){5.00}}
\put(25.00,70.00){\vector(-1,1){5.00}} \put(60.00,75.00){\circle*{2.00}}
\put(70.00,75.00){\circle*{2.00}} \put(65.00,75.00){\vector(-1,2){5.67}}
\put(59.33,86.00){\vector(-1,0){10.33}} \put(49.00,86.00){\vector(0,-1){22.00}}
\put(49.00,64.00){\vector(1,0){10.00}} \put(59.00,64.00){\vector(1,2){11.00}}
\put(70.00,86.00){\vector(1,0){11.00}} \put(70.00,64.00){\vector(-1,2){5.67}}
\put(81.00,86.00){\vector(0,-1){22.00}} \put(81.00,64.00){\vector(-1,0){11.00}}
\put(105.00,75.00){\circle*{2.00}} \put(115.00,75.00){\circle*{2.00}}
\put(109.00,75.00){\vector(0,1){11.00}}
\put(109.00,86.00){\vector(-1,0){15.00}}
\put(94.00,86.00){\vector(0,-1){22.00}} \put(94.00,64.00){\vector(1,0){17.00}}
\put(111.00,64.00){\vector(0,1){22.00}} \put(111.00,86.00){\vector(1,0){15.00}}
\put(126.00,86.00){\vector(0,-1){23.00}}
\put(126.00,63.00){\vector(-1,0){17.00}}
\put(109.00,63.00){\vector(0,1){12.00}} \put(15.00,45.00){\circle*{2.00}}
\put(25.00,45.00){\circle*{2.00}} \put(20.00,45.00){\oval(30.00,20.00)[]}
\put(39.00,45.00){\vector(1,0){7.00}} \put(65.00,45.00){\oval(30.00,20.00)[]}
\put(65.00,50.00){\circle*{2.00}} \put(65.00,40.00){\circle*{2.00}}
\put(65.00,40.00){\line(-1,1){5.00}} \put(65.00,50.00){\line(1,-1){5.00}}
\put(84.00,45.00){\vector(1,0){7.00}} \put(110.00,45.00){\oval(30.00,20.00)[]}
\put(15.00,16.00){\circle*{2.00}} \put(25.00,16.00){\circle*{2.00}}
\put(20.00,16.00){\oval(30.00,20.00)[]} \put(39.00,16.00){\vector(1,0){7.00}}
\put(65.00,16.00){\oval(30.00,20.00)[]} \put(65.00,11.00){\circle*{2.00}}
\put(84.00,16.00){\vector(1,0){7.00}} \put(110.00,16.00){\oval(30.00,20.00)[]}
\put(20.00,27.00){\vector(0,-1){22.00}} \put(65.00,21.00){\circle*{2.00}}
\put(65.00,27.00){\line(-1,-2){3.00}} \put(62.00,21.00){\line(3,-5){6.00}}
\put(68.00,11.00){\vector(-1,-2){3.00}} \put(105.00,16.00){\circle*{2.00}}
\put(115.00,16.00){\circle*{2.00}} \put(110.00,27.00){\line(0,-1){10.00}}
\put(110.00,17.00){\line(-5,2){5.00}} \put(105.00,13.00){\line(5,3){10.00}}
\put(115.00,13.00){\line(-5,2){5.00}} \put(110.00,15.00){\vector(0,-1){10.00}}
\put(137.00,104.00){\makebox(0,0)[cc]{{\large a)}}}
\put(137.00,74.00){\makebox(0,0)[cc]{{\large b)}}}
\put(137.00,44.00){\makebox(0,0)[cc]{{\large c)}}}
\put(137.00,15.00){\makebox(0,0)[cc]{{\large d)}}}
\put(105.00,16.00){\oval(6.00,6.00)[l]} \put(115.00,16.00){\oval(6.00,6.00)[r]}
\put(115.00,105.00){\circle*{2.00}} \put(105.00,105.00){\circle*{2.00}}
\put(70.00,105.00){\circle*{2.00}} \put(60.00,105.00){\circle*{2.00}}
\put(15.00,106.00){\vector(1,0){10.00}} \put(15.00,104.00){\vector(1,0){10.00}}
\put(60.00,106.00){\line(1,2){5.00}} \put(65.00,116.00){\vector(1,-2){5.00}}
\put(60.00,104.00){\line(1,-2){5.00}} \put(65.00,94.00){\vector(1,2){5.00}}
\put(105.00,104.00){\vector(-1,0){11.00}}
\put(94.00,104.00){\vector(0,-1){10.00}} \put(94.00,94.00){\vector(1,0){32.00}}
\put(126.00,94.00){\vector(0,1){10.00}}
\put(126.00,104.00){\vector(-1,0){11.00}}
\put(105.00,106.00){\vector(-1,0){11.00}}
\put(94.00,106.00){\vector(0,1){10.00}} \put(94.00,116.00){\vector(1,0){32.00}}
\put(126.00,116.00){\vector(0,-1){10.00}}
\put(126.00,106.00){\vector(-1,0){11.00}}
\put(14.00,45.00){\vector(-1,0){9.00}} \put(25.00,45.00){\vector(1,0){10.00}}
\put(60.00,45.00){\vector(-1,0){10.00}} \put(70.00,45.00){\vector(1,0){10.00}}
\put(117.00,45.33){\vector(1,0){8.00}} \put(102.00,44.67){\vector(-1,0){7.00}}
\end{picture}
\caption{Proofs of Lemmas 5 and 6} \label{defs}
\end{figure}
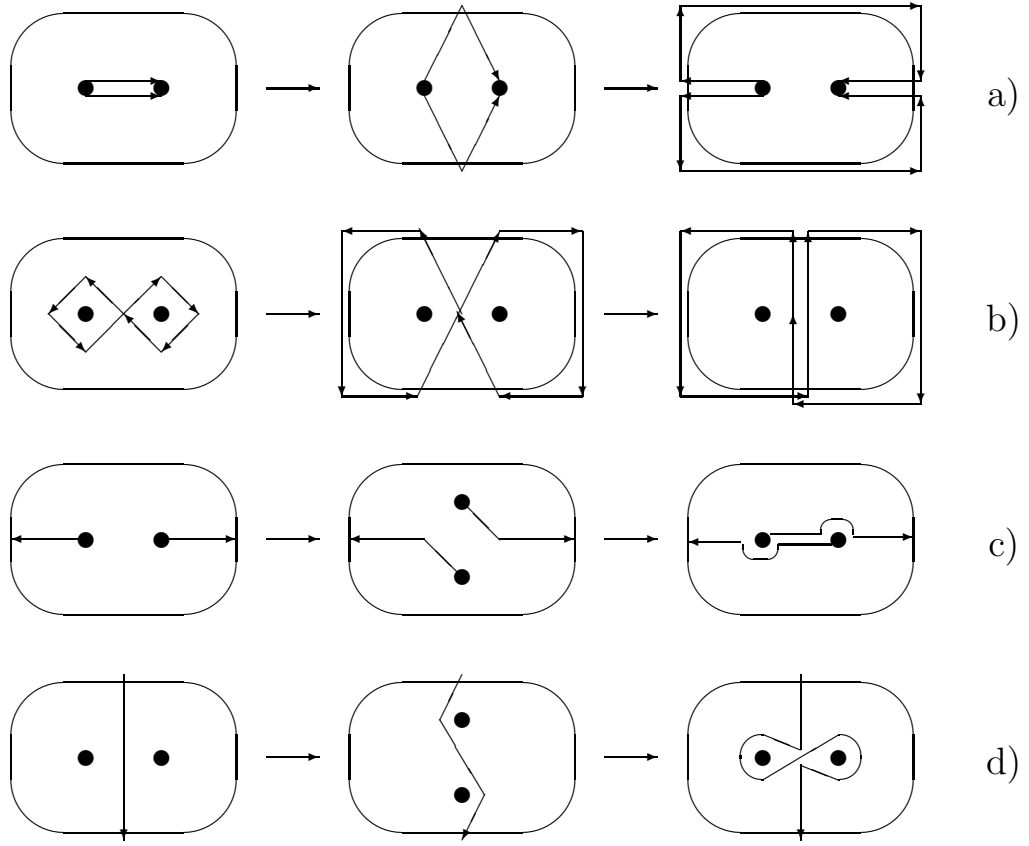

{\bf Lemma 5.} {\it For any $x \in {\mathcal H}_{i,\al}(r)$, the homomorphism $
\tilde J_{\al,r+1} \equiv \tilde j_{\al}\circ \tilde i_{\al}: {\mathcal
H}_{j+1,\al}(r+1) \to \bar {\mathcal H}_{j+1,\al}(r+1)$ (the reduction mod
$\partial B^{(r+1)}$) maps any element $\Sigma(x),$ $x \in {\mathcal
H}_{j,\al}(r),$ to $- \rlh (2x + \widetilde{Var}_{r} \circ \tilde J_{\al,r}(x))
\equiv -\rlh \circ (M_r+Id)(x)$.

The similar homomorphism $ J_{\al,r+1} \equiv j_{\al}\circ i_{\al}:
\chi_{j+1,\al}(r+1) \to \bar \chi_{j+1,\al}(r+1)$ maps any element $\infty
(y)$, $y \in \bar \chi_{j,\al}(r),$ to $- \dar (2y + J_{\al,r} \circ
Var_{r}(y)) \equiv - \dar \circ (\bar \mu_r + Id)(y)$. }
\medskip

{\it Proof.} See Figs. \ref{defs}a), b).
\medskip

{\bf Lemma 6.} {\it For any $r = m, \ldots, n-1,$ and any elements} $x \in
{\mathcal H}_{j,\al}(r)$ and $y \in \bar \chi_{j,\al}(r),$
$\widetilde{Var}_{r+1}(\rlh(x)) = \Sigma(x)$, $Var_{r+1} (\dar(y)) =
\infty(y)$.
\medskip

{\it Proof.} When $\xi$ moves along the circle $\delta \cdot e^{i\tau}, \tau
\in [0,2\pi],$ the ramification points $\pm \sqrt{\xi} \subset \Omega$ of the
fibre bundle $z_{n-r}: (B^{(r+1)},AX_\xi) \to \Omega$ move as is shown in Fig.
\ref{defs} c). We move the segments $[-\epsilon, -\sqrt{\xi}]$, $[\sqrt{\xi},
\epsilon]$ connecting the ramification points with ``infinity'' in such a way
that they coincide with the parts of the real axis close to the boundary of
$\Omega$. Let us also move the cycle realizing $\rightleftharpoons(\alpha)$ in
such a way that at any instant of this deformation it forms a fibre bundle with
standard fibre $\alpha/2$ over the union of two corresponding segments (except
for their endpoints). At the final instant of the monodromy each of these two
parts of $\rightleftharpoons(\alpha)$ will increase by the cycle
$\Sigma(\alpha)/2$ swept by the formal half of the initial cycle $\alpha$ in
transport over the added part of the resulting segment, so that $
\widetilde{Var}_{r+1}(\rightleftharpoons(\alpha)) = \Sigma(\alpha)$. For the
similar proof of the formula for the variation of $\dar(y)$, see Fig.
\ref{defs}d).
\medskip

This lemma proves statement b) of Theorem 3 and assertion of statement a)
concerning the variation of cycles $\rlh_l(x)$.
\medskip

Statement a) of Theorem 2 and assertion of statement b) concerning the action
on elements $\infty_l(y)$ follow from Lemmas 5 and 6 by induction over $r$.
\medskip

{\bf Corollary.} {\it The local monodromy operators $M_{r+1}, \bar M_{r+1},
\bar \mu_{r+1}$ and $\mu_{r+1}$ respectively map the elements $\Sigma(x),
\rlh(x), \dar(y)$ and $\infty(y)$ (where $x \in {\mathcal H}_{*,\alpha}(r), y
\in \bar \chi_{*,\alpha}(r)$) into $- \Sigma(M_r(x))$, $- \rlh(M_r(x))$, $-
\dar(\bar \mu_r(y))$ and $- \infty(\bar \mu_r(y))$ respectively.}
\medskip

{\bf Lemma 7.} {\it For any $x \in {\mathcal H}_{j,\al}(r)$ and} $y^* \in \bar
\chi_{2r-j,\al^*}(r)$,
$$\langle\Sigma(x), \dar(y^*)\rangle = (-1)^{1+j}\langle x, y^*\rangle ,$$
$$\langle\rlh(x), \infty(y^*)\rangle = (-1)^{1+j}\langle x, \bar
\mu_r(y^*)\rangle.$$

This lemma follows immediately from  constructions. Statements A) and B2) of
Theorem 4 follow from it by induction over $l$.
\medskip

{\bf Proposition 3.} {\it The homomorphisms $\Sigma : {\mathcal H}_{j,\al}(r)
\to {\mathcal H}_{j+1,\al}(r+1)$ and $\dar : \bar \chi_{j,\al}(r) \to \bar
\chi_{j+1,\al}(r+1)$ are isomorphisms for any $\al$ and $r=m, \ldots, n-1$.}
\medskip

{\it Proof.} Consider any smooth deformation contracting the disc $\Omega$ onto
the segment $[-\sqrt{\delta}, \sqrt{\delta}]$. Using the  fibre bundle
(\ref{przet}), we can lift this deformation to the space  $B^{(r+1)}\setminus
AX$. This lifted field allows us to realize any element of the group ${\mathcal
H}_{j+1,\al}(r+1)$ by a locally finite chain lying in this pre-image, and hence
also by one of the form $\Sigma(x),$ $x \in {\mathcal H}_{j,\al}(r)$. Thus the
map $\Sigma$ is epimorphic.

Our proposition follows now from the first equation of Lemma 7 and from the
fact that both Poincar\'e pairings ${\mathcal H}_{j,\al}(r) \otimes \bar
\chi_{2r-j,\al^*}(r) \to \C$ and ${\mathcal H}_{j+1,\al}(r+1) \otimes \bar
\chi_{2r-j+1,\al^*}(r+1) \to \C$ are non-degenerate.
\medskip

Theorem 1 is a direct corollary of this proposition.
\medskip

Denote by $\partial^- B^{(r+1)}$ the union of the usual boundary $\partial
B^{(r+1)}$ and the ``left half'' of $B^{(r+1)}$ consisting of points $z$ with
$z_{n-r} \le 0$.
\medskip

{\bf Proposition 4.} {\it For any $r=m, \ldots, n-1$ there is short exact
sequence}
\begin{equation}
\label{exact} 0 \to {\mathcal H}_{j}(r) \to \bar {\mathcal H}_{j+1}(r+1)
\stackrel{\rho}{\to} H^{lf}_{j+1}(B^{(r+1)}\setminus AX,
\partial^- B^{(r+1)}; L_\al) \to 0,
\end{equation}
{\it where the injection ${\mathcal H}_{j}(r) \to \bar {\mathcal H}_{j+1}(r+1)$
is the map $\rlh,$ and the epimorphism $\rho$ is induced by the reduction
modulo $\{z|z_{n-r} \le 0 \}$.}
\medskip

{\it Proof} (cf. \S \ 4.2 in \cite{V2}). We filter the disc $\Omega$ by the
sets $\{ \Psi_0 \subset \Psi_1 \subset \Psi_2 \}$ where $\Psi_0$ consists of
two points $\pm \sqrt{\delta}$, $\Psi_1$ consists of two segments
$[\sqrt{\delta}, \epsilon]$ and $[-\epsilon, -\sqrt{\delta}]$ (so that the set
$\Omega \setminus \Psi_1$ is a 2-cell), and $\Psi_2 \equiv \Omega$. Using the
projection (\ref{przet}) we lift this filtration onto the set
$B^{(r+1)}\setminus AX$; let $E_{p,q}^r$  be the spectral sequence generated by
this filtration and calculating the group $\bar {\mathcal H}_{*,\al}(r+1)$.
\medskip

{\bf Proposition 5.} {\it For any $q$, the elements $E_{p,q}^1$ of this
spectral sequence are as shown in the second row of Table 1.}
\medskip

\begin{table}
\caption{The groups $E^1_{p,q}$ for the main, anti-invariant and invariant
spectral sequences}
\begin{center}
\begin{tabular}{|c|cccc|}
\hline
$p$ & 0 & 1 & 2 & $\ge 3$ \\
\hline $E^1_{p,q}$ & $(H^{lf}_{q-1}(\partial {\mathcal F}_0),L_\al)^2$ & $(\bar
{\mathcal H}_{q,\al}(r))^2$ &
$\bar {\mathcal H}_{q,\al}(r)$ & 0  \\
\hline $E_{p,q}^{1,-}$ & $H^{lf}_{q-1}(\partial {\mathcal F}_0,L_\al)$ &
$\bar {\mathcal H}_{q,\al}(r)$ & 0 & 0 \\
\hline $E_{p,q}^{1,+}$  & $H^{lf}_{q-1}(\partial {\mathcal F}_0,L_\al)$ & $\bar
{\mathcal H}_{q,\al}(r)$ &
$\bar {\mathcal H}_{q,\al}(r)$ & 0 \\
\hline \end{tabular} \end{center} \end{table}

Indeed, the assertion concerning the group $E_{0,q}^1$ follows from the
statement b) of Lemma 2. For  any element $z \in \bar {\mathcal H}_*(r)$, two
corresponding elements in two summands $\bar {\mathcal H}_*(r)$ of $E_{1,*}^1$
are the homology classes of  two cycles (mod $z_{n-r}^{-1}(\Psi_0)$) swept out
by the copies of the cycle realizing $z$  transported over two segments  as in
the definition of the operation $\rlh$, see \S \ 3.3, and the corresponding
element of the group $E_{2,*}^1$ is swept out by the two-parametric family of
cycles obtained from the cycle realizing $z$ by the transportation to all
points of the 2-cell $\Omega \setminus \Psi_1$.
\medskip

Consider an involution acting on the term $E^1$ of our spectral sequence: it
permutes the elements of terms $E_{0,q}$, $E_{1,q}$, corresponding to the same
elements of groups $\tilde H^{lf}_{q-1}(\partial {\mathcal F}_0,L_\al)$, $\bar
{\mathcal H}_{q,\al}(r)$, and does not touch the terms $E_{2,q}$. Denote by
$E_{p,q}^{+,r}$ and $E_{p,q}^{-,r}$ the invariant and anti-invariant parts of
this involution, respectively.
\medskip

{\bf Lemma 8.} {\it The splitting of the term $E^1$ into the invariant and
anti-invariant parts is compatible with all further differentials and thus
defines the splitting of the entire spectral sequence. The groups $E^{+,1}$ and
$E^{-,1}$ of the invariant (respectively, anti-invariant) subsequence are as
shown in the fourth (respectively, third) row of the Table 1. The unique
nontrivial differential $\partial_1: \bar {\mathcal H}_{q,\al}(r) \to
H^{lf}_{q-1}(\partial {\mathcal F}_0,L_\al)$ of the anti-invariant subsequence
coincides with the boundary operator in ${\mathcal F}_0$, so that this
subsequence is nothing but the exact sequence of the pair $({\mathcal F}_0,
\partial {\mathcal F}_0)$ calculating the group $H^{lf}_*({\mathcal F}_0, L_\al) \equiv
{\mathcal H}_{q,\al}(r)$. The invariant subsequence is isomorphic to the
spectral sequence generated by the same filtration on the quotient space
$(B^{(r+1)}\setminus AX)/ (\partial^- B^{(r+1)}\setminus AX)$ and calculating
its homology group $H^{lf}_{i+1}(B^{(r+1)}\setminus AX,
\partial^- B^{(r+1)}\setminus AX; L_\al) $; this isomorphism is defined
by the reduction modulo $\{z|z_{n-r} \le 0\}$.}
\medskip

All this follows immediately from the construction and proves Proposition 4.
\medskip

By the construction, the induced involution of the limit homology group $\bar
{\mathcal H}_{i+1}(r+1)$ acts as multiplication by $-1$ on the subgroup
$\rlh({\mathcal H}_{i}(r))$ and acts trivially on the corresponding quotient
group, in particular the invariant subspace of this involution is canonically
isomorphic to $H^{lf}_{i+1}(B^{(r+1)}\setminus AX,
\partial^- B^{(r+1)}\setminus AX; L_\al) $.
\medskip

{\bf Lemma 9.} {\it Let $\beta$ be a sufficiently small closed disc in
$B^{(r+1)}$ centred at the distinguished point $\Delta^+$ (or $\Delta^-$) of
the fibre $\{z|z_{n-r} = \sqrt{\delta}\}$ (respectively, $-\sqrt{\delta}$).
Then the group $H^{lf}_{*}(B^{(r+1)}\setminus AX,
\partial^- B^{(r+1)}\setminus AX; L_\al) $
is isomorphic to $H_*^{lf}(\beta \setminus AX, \partial \beta; L_\al)$; this
isomorphism is induced by the reduction modulo the closure of the complement of
$\beta$.}
\medskip

This lemma follows from Lemma 2. Together with the K\"unneth decomposition
(\ref{kunn}) it proves the statement of Theorem 1$'$ concerning the group $\bar
{\mathcal H}_*$, and also the fact that the last summands in (\ref{stab2}) are
actually generated by the cycles described between Lemmas 4 and 5. Three
remaining assertions of Theorems 2--4 (triviality of  the action of
homomorphisms $J_{\al, (m+l)}$ (in Theorem 2) and $\widetilde{Var}_{(m+l)}$ (in
Theorem 3) on the last summands in (\ref{stab2}) and statement B1 in Theorem 4)
follow almost immediately from the construction of all these cycles. Lemma 4 is
just the fact Poincar\'e dual to the surjectivity of the map $\rho$ in
(\ref{exact}).

\end{document}